\documentclass[preprint,9pt]{article}

\usepackage{amssymb}
\usepackage{amsmath,booktabs,amsthm}
\usepackage{authblk}
\usepackage[dvipdfmx]{graphicx}

\newtheorem{thm}{Theorem}
\newtheorem{lem}{Lemma}
\newtheorem{cor}{Corollary}

\newtheorem{obs}{Observation}

\begin{document}

\title{A simple proof of tail--polynomial bounds on the diameter of polyhedra}

\author[1]{Shinji Mizuno\thanks{mizuno.s.ab@m.titech.ac.jp}}
\affil[1]{Department of Industrial Engineering and Management, Tokyo Institute of Technology, 
2-12-1-W9-58, Ookayama, Meguro-ku, Tokyo, 152-8552 Japan}

\author[2]{Noriyoshi Sukegawa\thanks{sukegawa@ise.chuo-u.ac.jp}}
\affil[2]{Department of Information and System Engineering, Chuo University, 
1-13-27 Kasuga, Bunkyo-ku, Tokyo 112-8551, Japan}

\maketitle

\begin{abstract}
Let $\Delta(d,n)$ denote the maximum diameter of a $d$-dimensional polyhedron with $n$ facets. 
In this paper, we propose a unified analysis of a recursive inequality about $\Delta(d,n)$ established by Kalai and Kleitman in 1992. 
This yields much simpler proofs of a tail--polynomial and tail--almost--linear bounds on $\Delta(d,n)$ which are recently discussed by Gallagher and Kim. 
\end{abstract}

\section{Introduction}\label{sec:intro}

The \emph{$1$-skeleton} of a polyhedron $P$ is an undirected graph $G=(V, E)$ which represents the vertex-vertex adjacency defined by the edges of $P$. 
More precisely, $V$ is a set of vertices of $P$, and $E$ is a set defined in such a way that $\{u, v \} \in E$ if and only if $\{ (1-\lambda) u +\lambda v: 0\le \lambda \le 1 \}$ forms an edge of $P$. 
The \emph{diameter} $\delta(P)$ of $P$ is the diameter of its $1$-skeleton. 
Formally, if we let $\rho_P(u, v)$ denote the shortest path length from $u$ to $v$, i.e., the number of edges required for joining $u$ and $v$ in $G$, then
$\delta(P)=\max\{ \rho_P (u, v) : u, v \in V \}$. 

Let $\Delta(d,n)$ denote the maximum diameter of a $d$-dimensional polyhedron with $n$ facets. 
In 1957, Hirsch conjectured that $\Delta(d,n) \le n-d$, which is disproved for unbounded polyhedra by Klee and Walkup~\cite{KlWa67}, and for even bounded polyhedra, i.e., polytopes, by Santos~\cite{Sa12}. 
An outstanding open problem in polyhedral combinatorics is to determine the behavior of $\Delta(d,n)$. 
In particular, the existence of a polynomial bound on $\Delta(d,n)$ is a major question. 
This question arose from its relation to the complexity of the simplex method; 
$\Delta(d,n)$ is a lower bound on the number of pivots required by the simplex method to solve a linear programming problem with $d$ variables and $n$ constraints. 

In this paper, we focus on the behavior of $\Delta(d,n)$ when $n$ is sufficiently large. 
Recently, in \cite{GaKi16}, Gallagher and Kim showed that 
\begin{thm}[Tail--polynomial bound~\cite{GaKi16}]\label{GaKi16}
For $d \ge 3$ and $n \ge 2^{d-1}$, 
\[\Delta(d, n) \le \frac{1}{16}\frac{n^3}{\sqrt{3\log(n)-5}}.  \]
\end{thm}\noindent
Their proof is based on two existing results: 
\begin{itemize}\setlength{\itemsep}{0pt}
\item[-] $\Delta(d, n) \le \Delta(d-1, n-1) + 2 \Delta\left(d, \left\lfloor \frac{n}{2} \right\rfloor \right) + 2$ for $\left\lfloor \frac{n}{2} \right\rfloor \ge d \ge 2$,
\item[-] $\Delta(d,n) \le \Delta(d,n+1)-1$ for $n>d>1$. 
\end{itemize}
The former is due to Kalai and Kleitman~\cite{KaKl92}, which has been used for proving bounds on $\Delta(d,n)$ for years, see, e.g., \cite{KaKl92,To14,SuKi15,Su16,GaKi16}. 
The latter is by Klee and Walkup~\cite{KlWa67}. 
It should be noted that Gallagher and Kim also show a tail--almost--linear bound $n^{1+\epsilon}$ for sufficiently large $n$. 

We propose a unified analysis of the former recursive inequality, Kalai--Kleitman inequality. 
As a corollary, we readily obtain tail--polynomial and tail--almost--linear bounds which are analogous to but slightly better than those discussed by Gallagher and Kim. 
Our proofs are much simpler and use Kalai--Kleitman inequality only. 

\section{Main results}

Using Kalai--Kleitman inequality, we show the following. 

\begin{thm}\label{Ours2}
Let $g(d,n)=\alpha^{d-3}(n-d)^k$ be a function with parameters $\alpha, k > 1$. 
If $\alpha$ and $k$ satisfy \[\frac{1}{\alpha}+\frac{1}{2^{k-1}} \le 1, \] then $\Delta(d,n) \le g(d,n)$ holds for $n \ge d \ge 3$. 
\begin{proof}
See Section~\ref{sec:Proof}. 
\end{proof}
\end{thm}\noindent
This bound $g(d,n)$ is motivated by Todd~\cite{To14}. 
He derived a bound $(n-d)^{\log_2(d)}$ by refining the proof of the sub-exponential bound $n^{\log_2(d)+2}$ shown by Kalai and Kleitman. 
As with Todd bound, our bound $g(d,n)$ takes 0 when $n=d$, and hence coincides with the fact that $\Delta(d,d)=0$. 

\begin{cor}\label{Cor2}
From Theorem~\ref{Ours2}, we obtain
\begin{itemize}
\item[\rm{a)}] for $d \ge 3$ and $n-d \ge 2^{d-3}$, $\Delta(d, n) \le (n-d)^3$, 
\item[\rm{b)}] for $d \ge 3$ and $n \ge 2^{2(d+2)}$, \[ \Delta(d, n) \le \frac{1}{16}\frac{n(n-d)^2}{\sqrt{3\log(n)-5}}, \]  
\item[\rm{c)}] for $\epsilon>0$, $d \ge 3$, and $n-d \ge \left[ \frac{2^{\epsilon/2}}{ 2^{\epsilon/2} -1} \right]^{2(d-3)/\epsilon}$, $\Delta(d, n) \le (n-d)^{1+\epsilon}$. 
\end{itemize}
\begin{proof}
Setting $\alpha, k:=2$ in Theorem~\ref{Ours2}, we have $\Delta(d,n) \le 2^{d-3}(n-d)^2$, which implies a). 
For b) and c), see Section~\ref{sec:Tail}. 
\end{proof}
\end{cor}\noindent
In particular, the bound in b) is analogous to but slightly better than that of Theorem~\ref{GaKi16}.

\section{Proofs of Theorem~\ref{Ours2}}\label{sec:Proof}

We always assume that $n \ge d$ since if otherwise the polyhedron has no vertex. 
Now, the goal is to show that 
\[ \text{P($d$):} \ \Delta(d,n) \le g(d,n)=\alpha^{d-3}(n-d)^k \ \text{for} \ n \ge d \]
is true for each $d \ge 3$. 
We prove this by induction on $d$. 
First, we can observe that P($3$) is true for any $\alpha, k > 1$. 
This is because we have $\Delta(3,n) \le n-3$, see, e.g., \cite{Kl64,Kl65,Kl66}. 

In what follows, assuming that P($d-1$) is true, we show that P($d$) is true. 
This is done by induction on $n$ while $d$ is fixed. 
When $n<2d$, it is known that $\Delta(d,n) \le \Delta(d-1,n-1)$, see, e.g., \cite{To14}. 
Hence, by the definition of $g$ and the validity of P($d-1$), we have \[ \Delta(d,n) \le \Delta(d-1,n-1) \le g(d-1,n-1) \le g(d,n) \]
for $n < 2d$. 
Then, let us consider the case when $n \ge 2d$. 
In this case, we employ the following result. 
\begin{lem}[Kalai--Kleitman inequality~\cite{KaKl92}]\label{thm:KaKl-ineq}
For $\left\lfloor \frac{n}{2} \right\rfloor \ge d \ge 2$, 
\[\Delta(d, n) \le \Delta(d-1, n-1) + 2 \Delta\left(d, \left\lfloor \frac{n}{2} \right\rfloor \right) + 2.  \]
\begin{proof}
See, e.g.,  \cite{KaKl92,To14}. 
\end{proof}
\end{lem}\noindent
In addition to the validity of P($d-1$), now, we can assume that P($d$) is true for $n'$ with $n' < n$ as the inductive hypothesis. 
Therefore, Lemma~\ref{thm:KaKl-ineq} implies 
\begin{align*}
\Delta(d,n)	&\le \Delta(d-1,n-1)+2\Delta\left(d, \left\lfloor \frac{n}{2} \right\rfloor \right)+2\\
				&= g(d-1,n-1)+2g\left(d, \left\lfloor \frac{n}{2} \right\rfloor \right)+2\\
				&\le \alpha^{d-4}(n-d)^k+2\alpha^{d-3}\left(\frac{n}{2}-d \right)^k +2\\
				&= \alpha^{d-3}(n-d)^k\left[ \frac{1}{\alpha}+\frac{1}{2^{k-1}}\left( 1 - \frac{d}{n-d}\right)^k+\frac{2}{\alpha^{d-3}(n-d)^k} \right] \\
				&= g(d,n)\left[ \frac{1}{\alpha}+\frac{1}{2^{k-1}} \left[ \left( 1 - \frac{d}{n-d} \right)^k+\left( \frac{1}{\alpha} \right)^{d-3}\left( \frac{2}{n-d}\right)^k \right] \right] \\
				&\le  g(d,n)\left[ \frac{1}{\alpha}+\frac{1}{2^{k-1}} \left[ 1 - \frac{d}{n-d} + \frac{2}{n-d} \right] \right] 
\end{align*}
where the last inequality follows from
\[ k>1,  \ \frac{1}{\alpha} <1, \ d \ge 3, \ \mbox{and} \ \frac{2}{n-d} \le 1. \]
Since $d \ge 3$,  we have 
\[ \Delta(d,n) \le g(d,n) \left[ \frac{1}{\alpha}+\frac{1}{2^{k-1}} \right] \le g(d,n). \]

\section{Comparison with Theorem~\ref{GaKi16} and a tail--almost--linear bound}\label{sec:Tail}

We add annotations to Corollary~\ref{Cor2}. 
First, let us observe b). 
Set $\alpha, k := 2$ as in a). 
Then, we have $\Delta(d,n) \le 2^{d-3}(n-d)^2$. 
Hence, it suffices to identify a lower bound $n_L$ such that $n \ge n_L$ implies 
\[ 2^{d-3} \le \frac{1}{16}\frac{n}{\sqrt{3\log(n)-5}}. \]
A sufficient condition for the inequality above is $2^{d+2} \le n/\sqrt{\log(n)}$. 
As $\sqrt{\log(n)} \le \sqrt{n}$ for $n \ge 1$, the condition can be further simplified to $2^{d+2} \le \sqrt{n}$. 
This means that we can set $n_L=2^{2(d+2)}$, which implies b). 

Then, let us observe c). 
For a given $\epsilon > 0$, set $k=1+\epsilon/2$. 
It is easy to see that $\alpha = \frac{2^{\epsilon/2}}{ 2^{\epsilon/2} -1}$ satisfies the condition $\frac{1}{\alpha}+\frac{1}{2^{k-1}} \le 1$. 
Hence, we have \[\Delta(d,n) \le \left[ \frac{2^{\epsilon/2}}{ 2^{\epsilon/2} -1} \right]^{d-3} (n-d)^{1+\epsilon/2}. \] 
Then, for $n$ satisfying $(n-d)^\frac{\epsilon}{2} \ge \left[ \frac{2^{\epsilon/2}}{ 2^{\epsilon/2} -1} \right]^{d-3}$, i.e., 
$n-d \ge \left[ \frac{2^{\epsilon/2}}{ 2^{\epsilon/2} -1} \right]^{\frac{2}{\epsilon}(d-3)}$, we have $\Delta(d,n) \le (n-d)^{1+\epsilon}$, which implies c). 

\section{Concluding remarks}

We finally point out that Larman bound~\cite{La70} $n2^{d-3}$ also derive tail--polynomial and tail--almost--linear bounds shown in Theorem~\ref{GaKi16}~\cite{GaKi16}. 
Note that however, Larman bound does not imply our bounds with the line of $n-d$. 
\begin{obs}\label{Obs1}
Since $\Delta(d,n) \le n2^{d-3}$ for $n \ge d \ge 3$~\cite{LaMaSa15}, we have 
\begin{itemize}\setlength{\itemsep}{0pt}
\item[-] for $d \ge 3$ and $n \ge 2^{d-3}$, $\Delta(d, n) \le n^2$, 
\item[-] for $\epsilon>0$, $d \ge 3$, and $n \ge 2^{\frac{1}{\epsilon}(d-3)}$, $\Delta(d, n) \le n^{1+\epsilon}$.  
\end{itemize}
\end{obs}\noindent
Larman bound was originally proven for only bounded polyhedra. 
Recently, in \cite{LaMaSa15}, Labb{\'e}, Manneville, and Santos proved it for simplicial complexes, which include general polyhedra. 



\begin{thebibliography}{9}

\bibitem{GaKi16}
J.M. Gallagher and E.D. Kim: 
Tail diameter upper bounds for polytopes and polyhedra. 
arXiv preprint arXiv:1603.04052. (2016). 
\bibitem{KaKl92}
G. Kalai and D.J. Kleitman: 
A quasi-polynomial bound for the diameter of graphs of polyhedra. 
\textit{Bulletin of the American Mathematical Society} 26 (1992) 315--316. 
\bibitem{Kl64}
V. Klee: 
Diameters of polyhedral graphs. 
\textit{Canadian Journal of Mathematics} 16 (1964) 602--614. 
\bibitem{Kl65}
V. Klee: 
Paths on polyhedra: I. 
\textit{Journal of the Society for Industrial and Applied Mathematics} 13 (1965) 946--956.
\bibitem{Kl66}
V. Klee: 
Paths on polyhedra: II. 
\textit{Pacific Journal of Mathematics} 17 (1966) 249--262.
\bibitem{KlWa67}
V. Klee and D.W. Walkup: 
The $d$-step conjecture for polyhedra of dimension $d < 6$. 
\textit{Acta Mathematica} 133 (1967) 53--78. 
\bibitem{La70}
D.G. Larman: 
Paths on polytopes. 
\textit{Proceedings of the London Mathematical Society} 20 (1970) 161--178. 
\bibitem{LaMaSa15}
J.-P. Labb{\'e}, T. Manneville, and F. Santos: 
Hirsch polytopes with exponentially long combinatorial segments. 
arXiv preprint arXiv:1510.07678. (2015).
\bibitem{Sa12}
F. Santos: 
A counter-example to the Hirsch Conjecture. 
\textit{Annals of Mathematics} 176 (2012) 383--412. 
\bibitem{SuKi15}
N. Sukegawa and T. Kitahara: 
A refinement of Todd's bound for the diameter of a polyhedron. 
\textit{Operations Research Letters} 43 (2015) 534--536. 
\bibitem{Su16}
N. Sukegawa: 
Improving bounds on the diameter of polyhedra in high dimensions. 
arXiv preprint arXiv:1604.04039. (2016). 
\bibitem{To14}
M.J. Todd: 
An improved Kalai--Kleitman bound for the diameter of a polyhedron. 
\textit{SIAM Journal on Discrete Mathematics} 28 (2014) 1944--1947. 
\end{thebibliography}
\end{document}